\documentclass{amsart}

\usepackage{amsmath,amssymb,amscd,amsxtra,amsthm}
\usepackage{latexsym}
\usepackage[dvips]{graphics,epsfig}

\newcommand{\R}{\mathbb{R}}

\newcommand{\nm}[1]{\|{#1}\|}

\newcommand{\ene}[3]{\mathcal{E}_{#1}(#2, #3)}

\newcommand{\V}[1]{Var_{\gamma}(|{#1}|^{2})}
\def\V#1{Var_{\gamma}({#1})}  

\newcommand{\Var}[1]{Var_{\gamma}(|{#1}|^{2})}
\def\Var#1{Var_{\gamma}({#1})}  
\newcommand{\dom}[1]{\mathcal{D}om{#1}}

\def\Si{Sierpi\'nski}
\def\Sig{Sierpi\'nski gasket}
\def\Sic{Sierpi\'nski carpet}
\def\pcf{p.c.f.}
\def\wup{weak uncertainty principle}
\def\erm{effective resistance metric}

\newtheorem{thm}{Theorem}

\newtheorem{prop}{Proposition}

\theoremstyle{remark}
\newtheorem{rem}{Remark}
\theoremstyle{definition}

\begin{document}

\title{Weak uncertainty principle for fractals, graphs and metric measure spaces}

\author{Kasso A.~Okoudjou}
\address{Kasso A.~Okoudjou\\
Department of Mathematics\\
University of Maryland\\
College Park, MD 20742-4015, USA} \email{kasso@math.umd.edu}

\author{Laurent Saloff-Coste$^\dag$}
\thanks{$^\dag$Supported in part by NSF grant DMS-0603886}
\address{Laurent Saloff-Coste\\
Department of Mathematics\\
 Malott Hall \\
Cornell
\newline University, Ithaca, NY 14853-4201, USA}
\email{lsc@math.cornell.edu}

\author{Alexander Teplyaev$^*$}
\thanks{$^*$Supported in part by NSF grant DMS-0505622}
\address{Alexander Teplyaev\\
Department of Mathematics\\
University of Connecticut\\
Storrs, CT 06269-3009, USA} \email{teplyaev@math.uconn.edu}


\date{\today}


\begin{abstract}
We  develop  a new approach to formulate and prove the weak
uncertainty inequality which was recently introduced by
Okoudjou and Strichartz.
We assume either an appropriate measure growth condition
with respect to the effective
resistance metric,  or, in the absence of such a metric, we assume the
 Poincar{\'e} inequality and reverse volume doubling  property.
We also consider the weak
uncertainty inequality in the context of Nash-type  inequalities.
Our results can be applied to a wide variety of metric
measure spaces, including graphs, fractals
 and manifolds.
\tableofcontents
\end{abstract}

\maketitle

\pagestyle{myheadings} \thispagestyle{plain} \markboth{K. A.
OKOUDJOU, L. SALOFF-COSTE, AND A. TEPLYAEV}{A WEAK UNCERTAINTY PRINCIPLE}

\section{Introduction}
The weak uncertainty inequality recently introduced in
\cite{OkSt04} for functions defined on \pcf\ fractals in general,
and on the \Sig\ in particular, obeys the same philosophy as the
classical uncertainty principle: it is impossible for any
normalized function to have a small energy and to be highly
localized in space. We refer to \cite{FoSi, PrRa85, Str89} for
more background on uncertainty principles. However, the existence
of localized eigenfunctions on some of these fractals (see
\cite{BaKi97, FuSh92, Sab00, Tep98}), is a main obstacle in
proving any analogue of the classical Heisenberg inequality. In
this paper we introduce a new approach to prove weak uncertainty
principles for functions defined on metric measure spaces equipped
with a Dirichlet (or energy) form,  which include certain fractals
and fractal graphs such as the Sierpi{\'n}ski lattice. More
precisely, we show that the weak uncertainty principle holds on
all spaces equipped with an effective resistance metric and a
measure satisfying an appropriate growth condition. Additionally,
we show that if instead of the existence of an effective
resistance metric on the space, we assume that a Poincar{\'e} or a
Nash  inequality holds along with another appropriate growth
condition on the measure, then it is also possible to prove the
weak uncertainty principle in this setting. In particular, our
results show that the self-similarity of the measure, which was
heavily used in \cite{OkSt04}, can be replaced by weaker
conditions. 

In order to formulate any  uncertainty inequality, one has to
 define  measures of space and frequency concentration.
For example, for complex-valued functions on $\mathbb R$ the
classical  Heisenberg Uncertainty Principle states that
$$
Var(|\hat{f}(\xi)|^{2})Var(|{f}(x)|^{2})\geqslant\frac{1}{16\pi^2} %
$$
for any function of $f \in L^2(\R)$ such that
$\nm{f}_{2}
=1$ and where $\hat{f}$ denotes the Fourier transform on $\R$.
This inequality can be rewritten in the   following form
$$
\int_{\mathbb R}\int_{\mathbb R}|x-y|^{2}|f(x)|^{2}\, |f(y)|^{2}\,
dx\, dy \int_{\mathbb R} |f'(x)|^{2}\, dx\geqslant\frac{1}{8}
$$
for any function of $L^2$ norm one. We refer to the survey article \cite{FoSi}
 for more information on the uncertainty principle.

In this paper we consider a metric measure space $(K, d, \mu)$, that is $(K, d)$  is a
metric space equipped with a Borel measure $\mu$. If
$\mathcal E$ is an energy form on this metric measure space, then we will
say that a weak uncertainty principle holds on $K$ if the
following estimate
\begin{equation}
\label{mainest-i} \V{u}\, \ene{}{u}{u} \geqslant C
\end{equation}
holds for any function $u \in L^{2}(K)\bigcap Dom(\mathcal E)$ such that  $\nm{u}_{L^2}=1$. Here $C$ is a
constant independent of $u$, and the spacial variance is defined
by
\begin{equation}\label{variance-i}
\V{u}=\iint_{K\times K}d(x, y)^{\gamma} |u(x)|^2\, |u(y)|^2\,
d\mu(x)\, d\mu(y).
\end{equation}

The central question of our paper is  the relation between
$d$, $\gamma$,
 $\mu$ and  $\mathcal E$ which implies the \wup, assuming that  the measure $\mu$
satisfies an appropriate growth condition. We formulate sufficient
conditions in several situations. The first one is when $d$ is the
so called  effective resistance metric on $K$ and satisfies
certain scaling properties. This setting is particularly relevant
in the context of analysis on fractals associated with
$\mathcal{E}$ and fractal graphs; see \cite{BB90,Ki2,Ki01,Ki03} for
more on the effective resistance metric. In this case
$\gamma=b+1$, where $b$ is an exponent which often  plays the role
of a dimension. Our result not only  provides a different and
simpler proof of \cite[Theorems 1 and 2]{OkSt04}, but also extends
them to all \pcf\ fractals \cite{Ki1,Ki01} and fractal graphs
\cite{BaBaKu1, GriTel01, GriHuLa03, HaKu1, KrTe03, MT1,MT2}, as
well as to modifications of them such as some fractafolds
\cite{Str99,Str03f}. Additionally, our result recovers the
classical  Heisenberg Uncertainty Principle in $\mathbb R$,
although not with the best constant. We also consider situations
where the effective resistance metric does not exist. In these
cases we assume some volume conditions and either that there is a
certain scaling in Poincar{\'e}'s inequality, or that a Nash-type
inequality holds. Either of these conditions allow us to prove our
result. These latter results are applicable for a wide variety of
metric measure spaces, ranging from graphs, to elliptic operators
on manifolds. Note that, in this case, the number $\gamma$
appearing in \eqref{variance-i} cannot, in general, be interpreted
as a dimension in the usual sense. However, $b=\gamma-1$ will
often represent the so-called walk dimension that appears
frequently in recent works on heat kernel estimates (see
\cite{BaBaKu1} and references therein). One of the features of our
results is their robustness. For example, since the \wup\ holds
for the \Si\ graphs, it also holds for the manifolds with similar
structures, e.g., the fractal-like manifolds considered in
\cite{KZ98}. Roughly speaking, the \wup\ holds if the energy and
measure have polynomial-type behavior, such as in the case of
fractals, fractal graphs and groups of polynomial growth.

Our paper is organized as follows. In Section~\ref{main} we state
and prove our main results. In particular, we prove that the weak
uncertainty principle holds under a variety of conditions raging
from the existence of an effective resistance metric, to assuming
Nash or Poincar\'e inequalities.
Section~\ref{examp} describes a few metric measure spaces for
which the main results of Section~\ref{main}
can be applied: \pcf\ fractals, uniform finitely ramified graphs,
\Sic s,\ fractal-like manifolds. We also discuss relation with
recent results on the heat kernel estimates on metric measure
spaces. 

Among the many other generalizations of the Heisenberg uncertainty
principle we wish to point out a recent preprint \cite{CRS}, in
which the authors work in the setting of Lie groups with
polynomial volume growth and which  partially overlap with our
results.

\medskip

\noindent {\bf Acknowledgments.} The authors are grateful to
Martin Barlow, Richard Bass  and Robert Strichartz  for many
helpful discussions, and to the anonymous referee(s) for
suggesting to consider Nash inequality in this context and
noticing the connection of our work with preprint \cite{CRS}.

\section{Main results}\label{main}

Let $(K, d)$ be a metric space equipped with a measure $\mu$ and  a
positive-definite symmetric quadratic form $\mathcal{E}$ with domain
$\mathcal{D}om(\mathcal{E}) \subset L^{2}(K)$.  Later, we will
impose some conditions relating the distance $d$ to the form
$\mathcal{E}$.
 We denote by  $B_r(x)$
 the ball with center $x$ and radius $r$ in the metric $d$.
To simplify notation, we make the convention that the $L^2$-norm
is infinite if a function is not square integrable, and that the
energy form is infinite if a function is not in its domain.

\subsection{Weak uncertainty principle and effective resistance metric}\label{spacerm}
In this subsection, we assume that $\mathcal{E}$ is a Dirichlet form and that the metric measure space $(K,
d, \mu)$ is such that $d$ is the effective resistance metric associated to  $\mathcal{E}$.
The most detailed study of such spaces, where $\mathcal{E}$ is also called a resistance form, is \cite{Ki03}.

The effective resistance metric  is defined by
$$
d(x,y)=\sup \mathcal E^{-1}(u,u),$$ where the supremum is
taken over all continuous functions $u$ such that $u(x)=1, u(y)=0$.  The
existence of the \erm\ is a nontrivial
 problem, see \cite{BB89,BB90,BB,BB2,Ki2,Ki01}; in particular, it is worth noticing
 that there are spaces without an effective resistance metric that is, for which the quantity above is infinite,
 e.g., for $\mathbb R^n$ with
$n\geqslant2$. However, on \pcf\ fractals and on some
Sierpi{\'n}ski carpets, which are ``not far from being one
dimensional'', it is known that the \erm\ does exist (see
Subsection~\ref{sics}).

We now state our first main result which generalizes Theorems $1$ and $2$ of \cite{OkSt04}.

\begin{thm}
\label{wupgene-p} Let $K$ be a space equipped with a measure $\mu$
and the effective resistance metric $d$ associated to the Dirichlet form $\mathcal{E}$. Assume that there exist positive
constants
$b, C_1, C_2$ such that  for all $x \in K$ and $r>0$
the following inequalities hold:
\begin{equation}
\label{measureequi} C_1r^b\leq \mu(B_r(x))\leq C_2r^b.
\end{equation}

 Then there exists $C>0$ such that for all $u
\in \mathcal{D}om(\mathcal{E}) $ with  $\nm{u}_{2}=1$ and $\gamma = b+1$ one has
$$\Var{u} \ene{}{u}{u} \geqslant C.$$
\end{thm}

\begin{proof}

Let $v=\Var{u}=\iint_{K\times K}d(x, y)^{\gamma} |u(x)|^2\,
|u(y)|^2\, d\mu(x)\, d\mu(y)$, where  $\gamma =
b+1$.
There exists $y$ such that
$$\int_{K} d^\gamma(x,y)u^2(x)d\mu(x) \leqslant v.$$
Let $r$ be defined by
\begin{equation}
\label{radius}
r=\sup\left\{s:\int_{B_s(y)}u^2(x)d\mu(x) < \frac12\right\}.
\end{equation}
For each  $s> 0 $ such that $\int_{B_{s}(y)}|u(x)|^{2}\, d\mu(x)
\leqslant \frac12$, we have
$$v\geqslant \int_{K - B_{s}(y)}d(x, y)^\gamma\, |u(x)|^{2}\, d\mu(x) \geqslant s^\gamma \,
\int_{K - B_{s}(y)}|u(x)|^{2}\,
d\mu(x)\geqslant \tfrac12s^\gamma,$$
 and the definition of $r$ implies that
\begin{equation}
\label{radvar}
r^\gamma\leqslant 2v.
\end{equation}
Moreover, by \eqref{measureequi} there is $c>1$ (it
suffices to take $c$ such that $c^{b} = 9 \tfrac{C_{2}}{C_{1}}$)
such that
\begin{equation}
\label{measbigbal}\mu(B_{cr}(y))>8C_2r^b.
\end{equation}

Let $t$ be an arbitrary number such that $t>r$. By the definition
of $r$ we have
\begin{equation}
\int_{B_t(y)}u^2(x)d\mu(x) \geqslant \frac12
\end{equation}
which together with \eqref{measureequi}
 yields
$$\sup_{B_t(y)}u^2(x)\geq\frac1{2C_2t^b}.$$
Additionally,  $\nm{u}_{2}=1$ yields
$$\inf_{B_{ct}(y)}u^2(x) < \frac1{8C_2t^b}.$$

Consequently,  there are $x_1,x_2\in B_{ct}(y) $ such that
$$(u(x_1)-u(x_2))^2\geq\frac1{8C_2t^b}.$$
Thus, by the definition of the  effective resistance metric,
$$\ene{}{u}{u}\geqslant \frac{(u(x_1)-u(x_2))^2}{d(x_1,x_2)}
\geq\frac1{16\, C_2\, c\, t\, t^b}= \frac1{16\, C_2\, c\,
t^{b+1}}.$$ Since this last inequality holds for all $t>r$ we
conclude that
$$  \ene{}{u}{u}\geqslant \frac1{16\, C_2\, c\, r^{b+1}}= \tfrac{C_{1}^{1/b}}{16\, C_{2}^{1 + 1/b}\,9^{1/b}_{\vphantom{f}}\, r^{b+1}}_{\vphantom{f}}.$$ Using
now \eqref{radvar} we conclude that
$$v\ene{K}{u}{u}
\geq \tfrac{C_{1}^{1/b} \, v}{16 \, C_{2}^{1 + 1/b}\, 9^{1/b}_{\vphantom{f}} \, v^{\frac{b+1}{\gamma}}}_{\vphantom{f}}=\tfrac{C_{1}^{1/b}\, v}{16\, C_{2}^{1 + 1/b}\, 9^{1/b}_{\vphantom{f}} \, v}=\tfrac{C_{1}^{1/b}}{16\, C_{2}^{1 + 1/b}\, 9^{1/b}_{\vphantom{f}} }=C,$$ where we have used the fact $\gamma=b+1$.
\end{proof}

\begin{rem}
 \noindent (a) The metric measure space $(K, d, \mu)$ in Theorem~\ref{wupgene-p} with the measure $\mu$ satisfying   ~\eqref{measureequi} is said to be an Ahlfors regular space \cite{Hei01}. Moreover, the constant $b$ appearing in
~\eqref{measureequi} is the Hausdorff dimension of $(K,
d)$.

\noindent (b)  Theorem~\ref{wupgene-p} holds when $K=\mathbb R^1$ with $b=1$ and is exactly the classical Heisenberg uncertainty principle of which we have given yet a different proof except for the precise value of $C$.

\noindent (c) One can modify this proof for various cases when  properties of the space are  different on  small and large scales. Two such modifications are given in
Theorems~\ref{wupgene-bdd} and~\ref{wupgene-grph}.
\end{rem}

Theorem~\ref{wupgene-p} assumes implicitly that the space $(K, d)$
is unbounded. In particular, it is not applicable to such
interesting examples as p.c.f. fractals. The next result which is
a weaker version of the previous one, generalizes the main result
of \cite{OkSt04} and deals with bounded spaces. We omit its proof
since it follows from obvious modifications from the proof of
Theorem~\ref{wupgene-p} below.

\begin{thm}
\label{wupgene-bdd} Let $K$ be a space equipped with a measure $\mu$
and the effective resistance metric $d$ associated to the Dirichlet form $\mathcal{E}$. Assume that there exist positive
constants
$b, C_{0}, C_1, C_2$ such that  ~\eqref{measureequi} holds for all $x \in K$ and all $0< r< C_0$.

Then there exists $C>0$  such that for all $u
\in \mathcal{D}om(\mathcal{E}) $ with  $\nm{u}_{2}=1$ and $\gamma = b+1$ one has
$$\Var{u} (\ene{}{u}{u} +1) \geqslant C.$$

\end{thm}

  Similarly,
Theorem~\ref{wupgene-p} excludes spaces where the local structure
is significantly different from the global one, for instance,
manifolds, graphs and spaces equipped with measure having
 atoms. In these cases the following variant of the our first
result can be proved using similar ideas.

\begin{thm}
\label{wupgene-grph}Let $K$ be a space equipped with a measure
$\mu$ and the effective resistance metric $d$ associated to the Dirichlet form $\mathcal{E}$. Assume that there
exist positive constants $b, C_{0}, C_1, C_2$ such that
~\eqref{measureequi} holds for all $x \in K$ and all $r> C_0$.
Then there exists  $C>0$  such that for
all $u \in \mathcal{D}om(\mathcal{E})$  with $\nm{u}_{2}=1$ one has
$$(\Var{u}+1) \ene{K}{u}{u} \geqslant C.$$
\end{thm}

\subsection{Weak uncertainty principle and Poincar\'e-type inequality}\label{spacepi}
In this subsection we no longer assume that $d$ is an effective
resistance metric. Instead, we will assume that a Poincar\'e
inequality holds. No specific property of the form $\mathcal{E}$ is required for the proof, it can be taken to be an
arbitrary positive-definite symmetric quadratic form on $L^{2}(K, d\mu)$. In this case the following theorem holds.

\begin{thm}
\label{wupgene} Let $K$ be a space equipped with a measure $\mu$
and a metric $d$ (not necessarily  an effective resistance
metric). Assume that there exists a positive constant $C_1$ such
that the energy form $\mathcal E$ on $K$
 satisfies the following Poincar{\'e} inequality for all locally square integrable functions
 $ u \in \mathcal{D}om(\mathcal{E})$
\begin{equation}
\label{poincare}
\int\limits_{B_r(y)}(u(x)-\bar{u}_{B_r(y)})^2d\mu(x) \leq
C_1r^{\gamma}\ene{}{u}{u},\end{equation} where $\bar{u}_{B_r(y)}$
is the average of $u$ over ${B_r(y)}$, and $\gamma $ is a positive
constant. Furthermore, assume that the measure $\mu$  satisfies
the ``reverse volume doubling property'' that is, there exist
an integer $k$ and a constant $C_2>1$ such that for all for all $x \in K$ and $r>0$
 we have
\begin{equation}
\label{inversedv} C_2\mu(B_r(x))\leq \mu(B_{kr}(x)).
\end{equation}

 Then there exists $C>0$ such that for all $u
\in \dom{\mathcal{E}}$ with $\nm{u}_{2}=1$   one has
$$\Var{u} \ene{}{u}{u} \geqslant C.$$
\end{thm}

\begin{proof}

 The proof is similar to that of Theorem~\ref{wupgene-p} with the main difference coming from the fact we need to
 use Poincar\'e's inequality instead of the basic properties of the effective resistance metric.

Let  $v=\Var{u}=\iint_{K\times K}d(x, y)^{\gamma} |u(x)|^2\,
|u(y)|^2\, d\mu(x)\, d\mu(y)$, and choose $y \in K$ such that
$\int_{K}d(x, y)^{\gamma}\, |u(x)|^{2}\, d\mu(x) \leq v$. Define  $r$ by
$$
r=\sup\left\{s:\int_{B_s(y)}u^2(x)d\mu(x) < \frac12\right\}.
$$ Note
that the definition of $r$ implies that $ r^\gamma\leqslant 2v.$
Now, iterating ~\eqref{inversedv} we can choose an integer $n>1$
and a constant $c$  such that  $ C_{2}^{n} \geq
16$ {and} $c \geq k^{n+1}$ where $C_2$ is the constant
appearing in ~\eqref{inversedv}. It suffices to choose $n = \lfloor
\tfrac{4\log 2}{\log C_{2}}\rfloor + 1$.
For this choice of $c$ we have,
\begin{equation}
\label{epsvolballs} T=\frac{\mu(B_{cr}(y))}{\mu(B_{kr}(y))} \geq
16.
\end{equation}
Notice that the definition of $r$  implies
that  $ \int\limits_{B_{kr}(y)}u^2(x)d\mu(x)\geq \tfrac{1}{2}$.
 Since $$ \left|\int\limits_{B_{cr}(y)}u(x)d\mu(x)\right| \leq \sqrt{\mu(B_{cr}(y))}$$ we conclude that
  $$|\bar{u}_{B_{cr}(y)}| \leq \tfrac{1}{\sqrt{\mu(B_{cr}(y))}}.$$
Now, using  \eqref{poincare} we see that
\begin{align*}
C_1c^\gamma r^\gamma\ene{}{u}{u} & \geq
\int\limits_{B_{cr}(y)}(u(x)-\bar{u}_{B_{cr}(y)})^2d\mu(x) \notag \\
&  = \nm{u-\bar{u}_{B_{cr}(y)}}_{L^2_{B_{cr}(y)}}^2\notag \\
& \geq \nm{u-\bar{u}_{B_{cr}(y)}}_{L^2_{B_{kr}(y)}}^2\notag \\
  &\geq \Big( \nm{u}_{L^{2}(B_{kr}(y))} -
\nm{\bar{u}_{B_{cr}(y)}}_{L^{2}(B_{kr}(y))}\Big)^{2}\\
& = \Big( \nm{u}_{L^{2}(B_{kr}(y))} -
|\bar{u}_{B_{cr}(y)}| \sqrt{\mu(B_{kr}(y))}\Big)^{2}\\
  & \geq \Big( \tfrac{1}{\sqrt{2}} - \tfrac{1}{4}\Big)^2 = \tfrac{9-4\sqrt2}{16}.
\label{fctoft}
\end{align*}
Thus, using the fact that $r^\gamma \leq 2 v$  we can write
$$v\ene{}{u}{u}
\geq\tfrac{9-4\sqrt2}{16}\frac v{C_1c^\gamma r^\gamma}\geq
\tfrac{9-4\sqrt2}{32}\frac v{C_1c^\gamma v} =
\tfrac{9-4\sqrt2}{32C_{1}c^\gamma}=C.$$ This last estimate
concludes the proof.
\end{proof}

It is interesting to observe that Theorem~\ref{wupgene} can be proved under a slightly different form of the
Poincar\'e inequality. In fact we have the following variant of Theorem~\ref{wupgene}.
Note that, for some applications,  Theorem~\ref{strongwupgene} is in a sense the strongest result of our paper,  as discussed in Subsection~\ref{subsection-Groups} in relation to the analysis on groups.

\begin{thm}\label{strongwupgene}
The conclusion of Theorem~\ref{wupgene} holds when condition ~\eqref{poincare}
is replaced by the following modified
Poincar\'e inequality
\begin{equation}\label{modpoin}
\int\limits_{B_r(y)}(u(x)-{u}_{r}(x))^2d\mu(x) \leq
C_{1}\,r^{\gamma}\ene{}{u}{u},\end{equation} where
$u_{r}(x)=\bar{u}_{B_{r}(x)}=\tfrac{1}{\mu(B_{r}(x))}\int_{B_{r}(x)}
u(y)\, d\mu(y),$ and $\gamma $ and $ C_1$ are some positive
constants.
\end{thm}

 \begin{proof} We define $v$,   $r$  and $y$ as in the proof  of Theorem~\ref{wupgene}. Then, again similarly to the proof  of Theorem~\ref{wupgene}, we can   iterate ~\eqref{inversedv} and find  a constant $c>2k$  such that
     $$    \frac{\mu(B_{cr/2}(y))}{\mu(B_{kr}(y))} \geq
     16.       $$
     Then for any  $z \in B_{kr}(y)$ we have
     $$|u_{cr}(z)| \leq \frac{1}{\sqrt{\mu\left(B_{cr}(z)\right)}}
 \leq \frac{1}{\sqrt{\mu\left(B_{\frac{cr}{2}}(y)\right)}}.$$
     Consequently, we have $$\nm{u_{cr}}_{L^{2}(B_{kr}(y))} \leq \sqrt{\frac{\mu(B_{kr}(y))}{\mu(B_{cr/2}(y))}} \leq \frac14.$$
     Now, using  \eqref{modpoin} we see that
     \begin{align*}
     C_1c^\gamma r^\gamma\ene{}{u}{u} & \geq
     \int\limits_{B_{cr}(y)}(u(x)-\bar{u}_{B_{cr}(y)})^2d\mu(x) \notag \\
       &\geq \Big( \nm{u}_{L^{2}(B_{kr}(y))} -
     \nm{u_{cr}}_{L^{2}(B_{kr}(y))}\Big)^{2}\\
       & \geq \Big( \tfrac{1}{\sqrt{2}} - \tfrac{1}{4}\Big)^2 = \tfrac{9-4\sqrt2}{16}.
     \end{align*} Then the proof then follows the same argument as in the proof of Theorem~\ref{wupgene}.
     \end{proof}

The following proposition gives a comparison between ~\eqref{poincare} and ~\eqref{modpoin} under the volume doubling condition.
\begin{prop}\label{compapoin}
Assume that the measure $\mu$ satisfies the doubling volume condition, that is, there exists $C>0$ such that for all $x\in K$ and all $r>0$ we have
\begin{equation}\label{vdc}
\mu(B_{2r}(x))\leq C \mu(B_{r}(x)).
\end{equation} Let $u$ be a locally square
integrable function that satisfies~\eqref{poincare}. Then $u$
satisfies also~\eqref{modpoin}.
\end{prop}

\begin{proof}
Let $u$ be a locally square integrable function  that satisfies~\eqref{poincare}. Observe that
$$
\|u-u_{r}\|_{L^{2}(B_{r}(y))} \leq \|u-\bar{u}_{B_{r}(y)}\|_{L^{2}(B_{r}(y))} +
\|u_{r}-\bar{u}_{B_{r}(y)}\|_{L^{2}(B_{r}(y))} =I +J.$$ By~\eqref{poincare} we see $I^{2} \leq C_{1} r^{\gamma}\ene{}{u}{u}, $  where $C_1$ is the constant appearing in ~\eqref{poincare}. To complete the proof it suffices to give a similar estimate for $J^2$. To simplify the notation, we let $V(x, r) = \mu(B_{r}(x))$ for $x \in K$ and $r>0$.
The second term in the right hand side of the last estimate can be estimated by
\begin{align*}
J^{2} & = \int_{B_{r}(y)}(u_{r}(x) - \bar{u}_{B_{r}(y)})^{2}\, d\mu(x)\\
& = \int_{B_{r}(y)} \Big| \tfrac{1}{V(x, r)}\int_{B_{r}(x)}u(t) d\mu(t) - \tfrac{1}{V(y, r)}
\int_{B_{r}(y)}u(z) d\mu(z)\Big|^{2}dx\\
&= \int_{B_{r}(y)} \Big|\tfrac{1}{V(x, r)V(y, r)}\iint_{B_{r}(x) B_{r}(y)}(u(t) - u(z) )d\mu(z)\, d\mu(t)\Big|^{2}dx.\\
& = \int_{B_{r}(y)}M(x) d\mu(x).
\end{align*}
For each $x \in B_{r}(y)$ we have by Jensen's inequality
\begin{align*}
M(x) &\leq \tfrac{1}{V(x, r)V(y, r)}\iint_{B_{r}(x) B_{r}(y)}|u(t) - u(z)|^{2}d\mu(z)\, d\mu(t)\\
& \leq\tfrac{1}{V(x, r)V(y, r)}\int_{B_{2r}(x)} \Big(\int_{B_{2r}(x)}|u(z) - u(t)|^{2}d\mu(z)\Big)d\mu(t) \\
& \leq\tfrac{4\,  V(x, 2r)}{V(x, r)V(y, r)}\, 2^{\gamma}\, C_{1}\, r^{\gamma}\, \ene{}{u}{u}\,
\end{align*}
where we have used once again ~\eqref{poincare}. By ~\eqref{vdc}, the last inequality becomes $$ M(x) \leq\tfrac{4}{V(y, r)}\, 2^{\gamma}\, C\, C_{1}\, r^{\gamma}\, \ene{}{u}{u}.$$
 Substituting this last inequality back in $J^2$ yields $$ J^{2} \leq \int_{B_{r}(y)}\tfrac{4 }{V(y, r)}\,\, 2^{\gamma}\, C\, C_{1} r ^{\gamma}\mathcal{E}(u, u) \, d\mu(x) = 4\, 2^{\gamma}\, C\, C_{1}\, r ^{\gamma}\, \ene{}{u}{u}.$$ which complete the proof of the proposition.
\end{proof}

 It is easily seen that both Theorems~\ref{wupgene} and~\ref{strongwupgene} only apply to unbounded spaces. For bounded spaces, we have the following modification of these results.

\begin{thm} \label{poinbound}
 Let $K$ be a space equipped with a measure $\mu$
and a metric $d$ (not necessarily  an effective resistance
metric). Assume that there exist positive constants $\gamma,
C_{0}, C_1$ and $ C_2 >1$ such that  ~\eqref{inversedv} together
with either ~\eqref{poincare} or ~\eqref{modpoin} hold for all $x
\in K$ and all $0< r< C_0$.

Then there exists $C>0$  such that for all $u
\in \mathcal{D}om(\mathcal{E}) $ with  $\nm{u}_{2}=1$ and $\gamma = b+1$ one has
$$\Var{u} (\ene{}{u}{u} +1) \geqslant C.$$
\end{thm}

\begin{rem}
\noindent (a) In general, and as opposed to the constant $b$
appearing in Theorem~\ref{wupgene-p}, the constant $\gamma$ in
\eqref{poincare} may not represent any sort of dimension
 of the space $(K, d)$ in the usual sense. However, in some of the examples we consider later,
 $\gamma$ is related to the so
 called
 walk dimension (see \cite{BB,BaBaKu1} and references therein).

 \noindent (b) It is worth pointing out that in all of the results we prove in this paper, the measure $\mu$
  is not necessarily a self-similar measure. But if it is self-similar, then the
measure weights and the resistance scaling weights are related by
a power law.

\noindent (c) Theorems~\ref{wupgene} and ~\ref{strongwupgene} apply to $K=\mathbb R^n $, $n\geqslant1$ with $\gamma=2$.

\noindent (d) One can show that on a metric measure space $(K, d,
\mu)$ where $d$ is the effective resistance metric, and where
$\mu$ satisfies ~\eqref{measureequi}, then the Poincar\'e
inequality ~\eqref{poincare} holds as well as ~\eqref{inversedv}.
 Indeed, if $u$ is a function which has
$L^2_\mu$-norm one in a ball of radius $r$ with respect to the
effective resistance metric, then \eqref{measureequi} implies that
the supremum of $|u|$ has a lower bound of the form $A \cdot
r^{-b/2}$ for some positive constant $A$. If, in addition, $u$ is
orthogonal to constants, then the difference between the supremum
of $u$ and  the infimum of $u$ also has a lower bound of the
 form $B \cdot r^{-b/2}$. Then the definition of the effective resistance implies that $\mathcal E(u,u)$ has
 a lower bound of the form $C \cdot r^{-b-1}$.

\noindent (e) If $(K, d)$ is a non-compact space and $d$ is
``geodesic'', that is, if for any two points $x, y \in K$ there is
a continuous curve $s: [0, 1] \to K$ such that $s(0) = x, \,
s(1)=y$ and $d(x, s(t)) = t$ then  ~\eqref{vdc} implies
~\eqref{inversedv}.
\end{rem}

\subsection{Weak uncertainty principle and Nash-type inequalities}\label{spacena}
In this subsection we investigate the relationship between
Nash-type inequalities and the weak uncertainty principle. In
particular, we first prove that the Nash inequality with small
enough dimension parameters implies the \wup\ if the measure is
upper $b$-regular. Furthermore, we use a Nash-type inequality to
prove a weak uncertainty principle even when these dimension
parameters are big. The Nash inequality plays an important role in
heat kernel estimates
\cite{BB3,BaBaKu1,Gri04,HaKu1,Ki01,Ki04,S-C02b} and references
therein.  Note that in general the Nash inequality does not imply
the Poincar{\'e} inequality.

\begin{thm}
\label{thmNash} Let $K$ be a space equipped with a  measure $\mu$,
a metric $d$ and a positive-definite symmetric quadratic form
$\mathcal E$ on $L^2_\mu$. Assume that the metric measure
space $(K, d, \mu)$ satisfies the following upper $b$-regularity
condition
\begin{equation}\label{nash-d-reg}
\mu(B_r(x))\leq C_1r^b, \, \textrm{for all}\,\, x \in K\, \textrm{and}\, \, r>0
\end{equation} and the following Nash
inequality
\begin{equation}\label{nash-nash}
\|f\|_{2}^{2+4/\theta}\leqslant
C_2\, \ene{}{f}{f}\, \|f\|_{1}^{4/\theta}
\end{equation}
for some positive constants $b, \theta, C_1$ and  $C_2$.
Furthermore assume that $\theta <2$. Then there exists $C>0$ such that for all $u \in \dom{\mathcal{E}} $ with
$\nm{u}_{2}=1$  and $\gamma = \tfrac{2b}{\theta}$ one has
$$\Var{u} \ene{}{u}{u} \geqslant C.$$ \end{thm}

\begin{proof}
Let $\nm{u}_{2}=1$ and  $v=\Var{u}=\iint_{K\times K}d(x,
y)^{\gamma} |u(x)|^2\, |u(y)|^2\, d\mu(x)\, d\mu(y)$. Then it suffices to prove that
\begin{equation}\label{eqnni}
\|u\|_1^{4/\theta}\leqslant C_{3} \cdot v,\end{equation} where $C_3$ is a constant to be specified.
By the definition of $v$, there
exists $y$ such that $\int_{K} d^\gamma(x,y)u^2(x)d\mu(x)
\leqslant v.$ For any $r>0$ we have

$$
\|u\|_1=\int_{B_r(y)}|u(x)|d\mu(x)+\int_{K-B_r(y)}|u(x)|d\mu(x)=A
+B.
$$
It is readily seen that $A \leq \big(\mu(B_r(y))\big)^{1/2} \leq
C_{1}^{1/2} r^{b/2}.$ Now we estimate $B$ as follows

\begin{align*}
B & = \int_{K-B_r(y)}|u(x)|d\mu(x)\\
 & = \int_{K-B_r(y)}d(x, y)^{-\gamma/2}\, d(x, y)^{\gamma/2}\,|u(x)|d\mu(x)\\
 & \leq \Big(\int_{K-B_r(y)}d(x, y)^{-\gamma}\, d\mu(x)\Big)^{1/2}\Big(\int_{K-B_r(y)}d(x, y)^{\gamma}\,|u(x)|^2
 d\mu(x)\Big)^{1/2}\\
&\leq \Big(\int_{K-B_r(y)}d(x, y)^{-\gamma}\,
d\mu(x)\Big)^{1/2}\Big(\int_{K}d(x, y)^{\gamma}\,|u(x)|^2
 d\mu(x)\Big)^{1/2}\\
 & \leq v^{1/2}\, \Big(\int_{K-B_r(y)}d(x, y)^{-\gamma}\,
 d\mu(x)\Big)^{1/2}.
 \end{align*}
 The integral in the last estimate can be estimated by

 \begin{align*}
 \int_{K-B_r(y)}d(x, y)^{-\gamma}\,
 d\mu(x)  & = \sum_{n \geq 0} \int_{\{ x: 2^{n} r \leq d(x, y) < 2^{n+1}
r\}}d(x, y)^{-\gamma}\, d\mu(x)\\
& \leq r^{-\gamma}\, \sum_{n \geq 0} 2^{-n \gamma}
\mu(B_{2^{n+1}r}(y))\\
& \leq C_{1}\, 2^{b}\,r^{-\gamma + b}\, \sum_{n \geq 0} 2^{n (b-\gamma)}\\
& = \tfrac{C_{1}2^{b}}{1 - 2^{b - \gamma}} r^{b-\gamma},
\end{align*}
where in the last estimate we have used the fact that $\theta < 2$ which is equivalent to $ b < \gamma$.  Therefore,

$$
\|u\|_1  \leqslant \tfrac{C_{1}^{1/2}2^{b/2}}{\sqrt{1 - 2^{b - \gamma}}} \, \big(r^{b/2}+ v^{1/2}r^{ b/2
-\gamma/2}\big).
$$
The minimum of the last expression with respect to $r$ is attained
when $$r=\big(\tfrac{\gamma -b}{b}\big)^{2/\gamma}v^{1/\gamma},$$
which implies
$$
\|u\|_1 \leqslant \tfrac{C_{1}^{1/2}2^{b/2}}{\sqrt{1 - 2^{b - \gamma}}}\, \tfrac {2\gamma -b}{\gamma -b}\,
\big(\tfrac{\gamma -b}{b}\big)^{b/\gamma}\,   v^{\tfrac{b}{2\gamma}},
$$
or equivalently $$ \|u\|_1 \leq C_{3}^{\theta/4}\, v ^{\theta/4}$$ where $C_{3}^{\theta/4} =
 \tfrac{C_{1}^{1/2}2^{b/2}}{\sqrt{1 - 2^{b - \gamma}}}\, \tfrac {2\gamma -b}{\gamma -b}\,
 \big(\tfrac{\gamma -b}{b}\big)^{b/\gamma},$ and $\frac{\theta}{4}= \frac{b}{2\gamma}$.
\end{proof}

 The following
variant of Theorem~\ref{thmNash} holds for all dimension
parameters.

\begin{thm}\label{monash}
Let $K$ be a space equipped with a  measure $\mu$, a metric $d$
and a symmetric quadratic  form $\mathcal E$ on $L^2_\mu$.
Assume that the metric measure space $(K, d, \mu)$ satisfies the
upper $b$-regularity condition ~\eqref{nash-d-reg},
 and the following
inequality
\begin{equation}\label{nashmod}
\|f\|_{2(1 + 2/\theta)}^{2+4/\theta}\leqslant
C_2\, \ene{}{f}{f}\, \|f\|_{2}^{4/\theta}
\end{equation}
for some positive constants $b, \theta, C_1$ and $C_2$. Then there
exists  $C>0$ such that for all $u \in \dom{\mathcal{E}} $ with $\nm{u}_{2}=1$  and $\gamma = \tfrac{2b}{\theta}$ one has
$$\Var{u} \ene{}{u}{u} \geqslant C.$$ \end{thm}

\begin{proof}
Let $\nm{u}_{2}=1$ and  $v=\Var{u}=\iint_{K\times K}d(x,
y)^{\gamma} |u(x)|^2\, |u(y)|^2\, d\mu(x)\, d\mu(y)$. We will
denote $p=1 + 2/\theta$. Choose $y \in K$ such that
$\int_{K} d^\gamma(x,y)u^2(x)d\mu(x) \leqslant v.$

Let $r$ be defined by
$$
r=\sup\left\{s:\int_{B_s(y)}u^2(x)d\mu(x) < \frac12\right\}.
$$
This
definition implies that $r^{\gamma} \leq 2v$. Let $t$ be an arbitrary number such that $t>r$. Then
$\int_{B_{t}(y)}|u(x)|^{2}\, d\mu(x) \geq \frac12.$ Consequently,
we see that $$\frac12 \leq \Bigg(\int_{K}|u(x)|^{2p}\,
d\mu(x)\Bigg)^{1/p} \Big(\mu\big(B_{t}(y)\big)\Big)^{1 - 1/p},$$ or equivalently, $$
\tfrac{1}{2^{p}} \leq \|u\|_{2p}^{2p}\, \Big(\mu\big(B_{t}(y)\big)\Big)^{p -
1}.$$ Using now the hypotheses of the theorem we obtain
\begin{align*}
\tfrac{1}{2^{p}} & \leq \|u\|_{2p}^{2p}\big(\mu(B_{t}(y))\big)^{p
-1}\\
& \leq 2\, C_2\, C_1^{p-1} \, \ene{}{u}{u} t^{b(p-1)}\\
& = 2\, C_2\, C_1^{2/\theta}\, \ene{}{u}{u} t^{2b/\theta}
\end{align*}
where $C_1$ is the constant appearing in
~\eqref{nash-d-reg}. Since this estimate holds for all $t>r$, and using the fact that $r^\gamma \leq 2v$, we
 conclude that

$$\tfrac{1}{2^{p}}  \leq 2\, C_2\, C_1^{2/\theta}\, \ene{}{u}{u} r^{2b/\theta}
\leq 2^{1 + 2b/\theta \gamma}\, C_2\, C_1^{2/\theta}\,
\ene{}{u}{u}\, v^{2b /\theta\gamma}.$$  Therefore, $$\ene{}{u}{u}\, \Var{u}^{2b
/\theta\gamma} \geq 2^{-p -1 -2b/\theta \gamma }
C_{2}^{-1}C_{1}^{-2/\theta}=C.$$ This last inequality completes the proof by the choice of $\gamma$.
\end{proof}

\begin{rem} We wish to recall that in Theorems~\ref{thmNash} and~\ref{monash} the form $\mathcal{E}$ is not assumed
to be a Dirichlet form.
However, if $\mathcal{E}$ is a Dirichlet form, then it is known that  ~\eqref{nash-nash} is equivalent to
 ~\eqref{nashmod} for all $\theta >0$; see \cite{bclsc}. Moreover, in this case, it follows from
 \cite[Theorem 3.1.5]{S-C02b} that we can estimate the volume of a ball from below, i.e., for all $x \in K$
 and $r>0$ we
have $\mu(B_{r}(x)) \geq C r^{\theta}$, where $C$ is a positive
constant. Thus for all $r>0$ we have $C r^{\theta}
\leq \mu(B_{r}(x)) \leq C' r^{b}$, and so $\theta = b$. Consequently, if we assume in Theorem~\ref{thmNash} that
the energy form $\mathcal{E}$ is a Dirichlet form, then we can remove the restriction $\theta < 2$ by using the above
observation. Note that in this case, $\gamma =2$.
\end{rem}

\subsection{Local weak uncertainty principle}\label{local}

In this subsection we state some local versions of the results
proved above. The proofs are easy adaptation of those given above
and are omitted.

\begin{thm}
\label{locthmNash} Let $K$ be a space equipped with a  measure
$\mu$, a metric $d$ and a positive-definite symmetric quadratic
form $\mathcal E$ on $L^2_\mu$. Assume that the metric measure
space $(K, d, \mu)$ satisfies the following upper $b$-regularity
condition ~\eqref{nash-d-reg}  and the following Nash inequality
\begin{equation}\label{locnash-nash}
\|f\|_{2}^{2+4/\theta}\leqslant
C_2\big(\ene{}{f}{f}+\|f\|_2^2\big)\|f\|_{1}^{4/\theta}
\end{equation}
for some positive constants  $b, \theta, C_1$, and $ C_2$ with
$\theta <2$. Then there exists positive constants $C>0$ such that
for all $u \in \dom{\mathcal{E}}$ with $\nm{u}_{2}=1$ and $\gamma
= \tfrac{2b}{\theta}$ one has
$$\Var{u} (\ene{}{u}{u}+1) \geqslant C.$$ \end{thm}

Similarly, we have

\begin{thm}\label{locmonash}
Let $K$ be a space equipped with a  measure $\mu$, a metric $d$
and a symmetric quadratic  form $\mathcal E$ on $L^2_\mu$. Assume
that the metric measure space $(K, d, \mu)$ satisfies the upper
$b$-regularity condition ~\eqref{nash-d-reg},
 and the following
inequality
\begin{equation}\label{locnashmod}
\|f\|_{2(1 + 2/\theta)}^{2+4/\theta}\leqslant
C_2\big(\ene{}{f}{f}+\|f\|_2^2\big)\|f\|_{2}^{4/\theta}
\end{equation}
for some positive constant $b, C_1$ and $C_2$. Then there exists $C>0$ such that for all $u \in \dom{\mathcal{E}}$ with $\nm{u}_{2}=1$  and $\gamma = \tfrac{2b}{\theta}$ one has
$$\Var{u} (\ene{}{u}{u}  + 1)\geqslant C.$$ \end{thm}

\section{Applications and examples}\label{examp}

\subsection{\Sig\ and p.c.f. fractals}
As mentioned in the Introduction, the weak  uncertainty principle
for functions defined on the \Sig\, was first introduced in
\cite{OkSt04}. While the results in that paper were stated for
p.c.f. fractals, they were only proved for the \Sig. In this
subsection, we use the results of Section~\ref{main} not only to
provide a simpler proof to the main results of \cite{OkSt04}, but
also to establish weak uncertainty principles on all \pcf\
fractals. We briefly define the \Sig\ which is a typical
example of a p.c.f.\ fractal, and refer to \cite{Ba98,Ki01, Str99, Str05} for more background on analysis on p.c.f.\ fractals.

Consider the contractions maps $F_1, F_2$ and $F_3$ defined on
$\R^2$ by $F_1(x) = \frac{1}{2}x$, $F_2(x) = \frac{1}{2}x +
(\frac{1}{2}, 0)$ and $F_3(x) = \frac{1}{2}x + (\frac{1}{4},
\frac{\sqrt{3}}{4})$, for $x \in \R^2$. The Sierpi{\'n}ski gasket
$K=SG$, is the unique nonempty  compact
 subset of $\R^2$ such that
\begin{equation}\label{ifs}
K= \bigcup_{i=1}^{3}F_{i}K.
\end{equation}
For any positive integer $m$,
 $\omega=(\omega_{1}, \omega_{2}, \hdots, \omega_{m})$ where each $\omega_{i} \in \{1, 2, 3\}$ is called a word of
 length $|\omega|=m$, and we denote
 $F_{\omega}=F_{\omega_{m}}\circ F_{\omega_{m-1}}\circ \hdots \circ F_{\omega_{1}}$.
Then $F_{\omega}K$ is called a cell of level $m$ if $\omega$ is
 a word of length $m$. It is worth mentioning that $SG$ can be defined as a limit of graphs: Let $\Gamma_0$ be the
 complete graph with vertices $V_0=\{(0, 0), (1, 0), (\frac{1}{2}, \frac{\sqrt{3}}{2})\}$ which  are the fixed points
 of the contractions $F_i$. The graph $\Gamma_m$ with vertices $V_m$ is defined inductively by
$V_{m}= \bigcup_{i=1}^{3} F_{i}V_{m-1}$, $m\geqslant 1$, and $x
\sim_{m} y$ if $x$ and $y$ are in the same $m$-cell. The
(standard) measure on $K$ is the probability measure on $K$ that
assigns
 to each cell of level $m$ the measure $3^{-m}$.
It follows that $SG$ is equipped with a self-similar measure that
satisfies trivially \eqref{measureequi}. By  defining an energy
form on $SG$,  it can be shown that this gives rise to a
resistance metric on $SG$, see \cite{Ba98, Ki1, Ki01}.
Consequently,  Theorem~\ref{wupgene-bdd} recovers one of the main results -- Theorem $1$ -- of \cite{OkSt04} for $SG$.

More generally, let $\{F_{i}\}_{i=1}^{N}$ be a set of contractive injective
maps on a compact metric space $K$ and $$K=\bigcup_{i=1}^{N}F_i(K).$$ We assume that $K$ is a \pcf\ self-similar set in the sense of \cite{Ki1, Ki01}. Following \cite{Ki1, Ki01}, one can sometimes define a
self-similar  energy form
$\ene{}{\cdot}{\cdot}$ such that $$ \ene{}{u}{u}= \sum_{i=1}^{N}
\rho_{i} \ene{}{u\circ F_i}{u \circ F_i},$$ which gives
rise to an effective resistance metric. If the resistance scaling factors $\{\rho_i\}_{i=1}^{N}$ satisfy
the regularity condition $$\rho_i>1$$ then the  effective resistance metric induces the same topology as the original metric on $K$. Note that the energy form $\ene{}{\cdot}{\cdot}$ is defined without reference to any measure on $K$. If we also
consider a set of positive real numbers $\{\mu_{i}\}_{i=1}^{N}$, called the measure scaling factors, such that $\sum_{i=1}^{N}\mu_{i} = 1$, then we have a self-similar
measure $\mu$ on $K$ such that $\mu = \sum_{i=1}^{N} \mu_{i} \mu \circ
F_{i}^{-1}$. According to \cite{Ki01,Ki03}, the  dimension  of $K$ in the effective resistance metric is the unique $b$ such that $$\sum_{i=1}^{N}\rho^b_{i} = 1.$$ Therefore, the most natural choice of the measure scaling factors is $$\mu_i=\rho^b_{i}.$$ Notice that taking $\mu_i=\rho_{i}^{b}$ is also natural because with this choice the asymptotic behavior of the Weyl function is
well studied \cite{KL,Ki01}. Our Theorem~\ref{wupgene-bdd} holds on $K$ if the condition \eqref{measureequi} is satisfied.
The difficulty is, however, that if the measure scaling factors are not equal and some combinatorial-geometric conditions are not satisfied, then in many cases the self-similar measure $\mu$ does not satisfy even the volume doubling property much less the regularity condition \eqref{measureequi}. Fortunately, for many cases the  the regularity condition \eqref{measureequi} does hold, such as in the situation of the nested fractals. Many related question, and the Nash inequality in this context in particular, are discussed in \cite{Ki04}.

To construct non compact fractals which satisfy conditions of Theorem~\ref{wupgene-p}, we assume for simplicity that $\{F_{i}\}_{i=1}^{N}$ are contractive injective
maps on $\R^d$ and $K$ is the unique compact set such that $K=\bigcup_{i=1}^{N}F_i(K) $.
Then  one can construct an
increasing  sequence of sets $K_n$ using inverse maps  $F_{i_n}^{-1}$, and define the blowup of $K$
to be $K_{\infty}= \bigcup_{n=0}^{\infty}K_n$ where $K_0=K$. Then
$K_{\infty}$ is an unbounded self-similar set, called fractal
blowup and was first introduced in \cite{Str96}, see also
\cite{Sab00, Sab03, Sab04b, Tep98} for more about fractal blowups.

\subsection{\Si\ graphs}
As another application of the results of Section 2, we prove a
weak uncertainty principle of some graphs related to the \Sig\
$K=SG$, and its blowup $K_{\infty}$. More precisely, for any
integer $m \geqslant 0$, let $\Gamma_{m}$ be the $mth$ pre-gasket
approximation to $K$, i.e., the $mth$ graph approximation of $K$.
We define a (finite) graph $\Gamma_{-m}$ by
$\Gamma_{-m}=F_{\omega_{i}}^{-1} \circ F_{\omega_{2}}^{-1}\hdots
\circ F_{\omega_{m}}^{-1}(\Gamma_{m})$, and an infinite graph
$\Gamma_{\infty}$ by $\Gamma_{\infty}=\bigcup_{m\geqslant
0}\Gamma_{-m}$. The graph $\Gamma_{\infty}$ is an example of an infinite
self-similar graph, which is also referred to as the
Sierpi{\'n}ski lattice; we refer to \cite{Tep98} and the
references therein for more on this type of graphs. Note that for
all integers $m\geqslant 0$, $\Gamma_{-m}$ is similar to the
(finite) graph obtained by taking $F_{\omega_{i}}=F_1$ for all
$i$, in which case, $\Gamma_{-m}=2^{m}\Gamma_{m}$. In this case it
can be shown that  Theorem \ref{wupgene-grph} holds on
$\Gamma_{\infty}$.

\subsection{Uniform finitely ramified  fractals and graphs}
The uniform finitely ramified fractals (u.f.r.\ ) and the
unbounded fractals associated to them were introduced in
\cite{HaKu1} and include the nested fractals and are contained in
the class of \pcf\  self-similar sets, see \cite{HaKu1, Ki01,
Ki1}. Clearly, Theorems~\ref{wupgene-p} and~\ref{wupgene-bdd}
applies to these class of fractals.

Additionally, the (infinite) u.f.r.\ graphs were constructed from u.f.r.\ in
\cite{HaKu1}, where it was proved  that there exists an effective resistance metric on
this class of graphs. Therefore, using \cite[Lemma 3.2]{HaKu1}
 one can show that Theorem~\ref{wupgene-grph} applies in this
setting as well.

\subsection{\Sic s and graphical \Sic s}\label{sics}
These are examples of non finitely ramified fractals and fractal
graphs \cite{BB89,BB90,BB,BB2,KZ92}. In particular, they are non
p.c.f.\ fractals, and it is interesting to notice that most of our
results apply  in this setting. Hence, we answer affirmatively a
question posed in \cite{OkSt04} of whether the main results of
that paper apply to ``genuine'' non-p.c.f. fractals. More
precisely, on  the generalized Sierpi\'nski carpets (GSC) and the
unbounded sets that can be constructed based on them, it is known
that a two sided heat kernel estimate holds, \cite{BB89, BB,
BB92}. Thus, following \cite{BB89, BB, BB92,BaBaKu1} or
\cite[Theorem 3.2]{GriHuLa03}, one can show that
~\eqref{measureequi} holds on the GSC and all related sets; this
in turn implies that  \eqref{inversedv} holds also in these
settings. Moreover, it is known that the two sided heat kernel
estimate implies that the Poincar\'e inequality holds e.g., see
\cite{BB2, BaBaKu1}. Consequently, Theorem~\ref{wupgene} applies to all unbounded spaces
constructed on GSC, while Theorem~\ref{poinbound} applies to all Sierpi\'nski carpets for which such estimates exist. Moreover, Theorem~\ref{wupgene-p} applies to resistance Dirichlet forms on the  Sierpi\'nski carpets in dimension less than 2, such as self-similar Dirichlet forms on the  Sierpi\'nski carpets constructed in \cite{KZ92}.


\subsection{Groups.}\label{subsection-Groups}

Let us briefly describe how our results apply in the case when the underlying space $K=G$ is a group and the
distance $d$ and quadratic form $\mathcal{E}$ have some invariance properties.

For instance, let $G$ be a real connected Lie group equipped with
a left-invariant Riemannian structure (in fact, a left-invariant
sub-Riemannian structure would work as well); see, e.g.,
\cite{S-C02b, vscc} for details. We let $d$ be the Riemannian
distance and $\ene{}{f}{f}= \int_{G}|\bigtriangledown f|^{2}\,
d\mu$ where $\mu$ is the (left-invariant) Riemannian measure and
$|\bigtriangledown f|^{2}$ the Riemannian length of the gradient
of $f$. In this setting, Theorem~\ref{wupgene-p} and
Theorem~\ref{thmNash} apply only to the case $G=\R$ because in
higher dimensions the resistance metric is infinite.
 Observe that the reverse doubling condition ~\eqref{inversedv} holds on any non-compact Lie group. Therefore,
Theorem~\ref{wupgene} applies with $\gamma =2$ to all the cases where $G$ is non-compact with polynomial volume growth.
This includes $\R^n$ and all nilpotent Lie groups, in particular the Heisenberg group
$$\Bigg\{\left(\begin{array}{ccc}1& x& z\\ 0& 1 & y\\ 0& 0& 0\end{array}\right): x, y, z \in \R\Bigg\}.$$
We do not know if condition ~\eqref{poincare} implies polynomial growth, but it seems unlikely that it holds for
group of exponential growth.

 Theorem~\ref{strongwupgene} has the great advantage over Theorem~\ref{wupgene}
that it applies to all non-compact unimodular Lie groups. This because ~\eqref{modpoin} holds with $\gamma=2$ on any
such Lie group, see e.g., \cite[Theorem 3.3.6]{S-C02b}.

Finally, because of Remark $3$, Theorem~\ref{monash} is
essentially restricted to the case when $G=\R^n$ for some
$n=\theta$.

Consequently, in the case of Lie groups,
Theorem~\ref{strongwupgene} is by far the most powerful result.
Moreover, Theorem~\ref{strongwupgene} applies as well to the case
of infinite, finitely generated groups equipped with the counting
measure and a quadratic form of the type
$$\ene{}{f}{f}= \tfrac{1}{2|S|}\sum_{x \in G}\sum_{s \in S}|f(xs) - f(x)|^{2}$$ where $S=S^{-1}$ is a finitely
symmetric generating set. Indeed, in this setting, ~\eqref{modpoin} holds with $\gamma=2$
(the proof of \cite[Theorem 3.3.6]{S-C02b} may easily be adapted to this setting).

\subsection{Metric measure spaces and heat kernel estimates}
Our results of Section~\ref{main} are applicable to the general
setting of metric measure spaces. For a metric measure space $(K,
d, \mu)$ the main assumption we make  is the existence of a heat
kernel $\{p_{t}\}_{t > 0}$, which is the fundamental solution  of
the heat equation where the self-adjoint operator associated with
the energy form $\mathcal E$ plays the role of a Laplacian. If the
heat kernel, which is a  non-negative measurable function
$p_{t}(x, y)$ on $[0,\infty)\times K\times K$,  satisfies the
following two sided estimate for $\mu$-almost $x, y \in K$ and all
$t \in (0, \infty)$:
\begin{equation}\label{heatkest}
\frac{1}{t^{\alpha /\beta}}\Phi_{1}\bigg(\frac{d(x,
y)}{t^{1/\beta}}\bigg) \leqslant p_{t}(x, y) \leqslant \frac{1}{t^{\alpha
/\beta}}\Phi_{2}\bigg (\frac{d(x, y)}{t^{1/\beta}}\bigg),
\end{equation} where $\alpha$ is the Hausdorff dimension of $(K,
d)$ and $\beta = \alpha + 1$, and $\Phi_{1}, \Phi_{2}$ are
non-negative monotone decreasing functions on $[0, \infty)$, then
under a mild decay condition on $\Phi_2$, it is shown in
\cite[Theorem 3.2]{GriHuLa03} that \eqref{heatkest} implies
~\eqref{measureequi} with $b=\alpha$. This can be used in turn to
prove \eqref{inversedv}. Consequently, once a Poincar{\'e}-type
estimate is established in this setting, our results can be
applied.

Moreover, heat kernel estimates of the type
\begin{multline*}
\frac{c_1}{
\mu(B(x, t^{1/\gamma}))}
\exp {\left(-\left(  \frac{d(x,y)^\gamma}{c_1t}\right)^{\frac1{\gamma-1}}\right)}
\leqslant
p_{t}(x, y)\\
\leqslant
\frac{c_2}{
\mu(B(x, t^{1/\gamma}))}
\exp {\left(-\left(  \frac{d(x,y)^\gamma}{c_2t}\right)^{\frac1{\gamma-1}}\right)}
\end{multline*}
imply the
Poincar{\'e} inequality, and these estimates can be established on many fractals and other spaces (see \cite{BaBaKu1}
and references therein).


\begin{thebibliography}{13}

\bibitem{bclsc}
D.~Bakry, T.~Coulhon, M.~Ledoux, and L.~Saloff-Coste, {\it Sobolev
inequalities disguise,} Indiana Univ.\ Math.\ J.\ {\bf 44} (1995),
no.\ 4, 1033-1073.

\bibitem{Ba98} M.~T.~Barlow, Diffusion on Fractals, in: Lectures Notes in Mathematics, Vol.~1690, Springer,
Berlin, 1998.

\bibitem{BB89} M. T. Barlow and R. F. Bass,
\emph{The construction of Brownian motion on the \Sic.}\
Ann. Inst. H. Poincar\'{e} Probab. Statist.  \textbf{25} (1989), 225--257.

\bibitem{BB90} M. T. Barlow and R. F. Bass,
\emph{On the resistance of the \Sic.}\
Proc. Roy. Soc. London Ser. A   \textbf{431} (1990), 345--360.

\bibitem{BB} M. T. Barlow and R. F. Bass,
\emph{Brownian motion and harmonic analysis on Sierpinski
carpets.}\
Canad. J. Math., \textbf{51} (1999), 673--744.

\bibitem{BB2} M. T. Barlow and R. F. Bass,
\emph{Random walks on graphical {S}ierpinski carpets,}. Random walks and discrete potential theory
(Cortona, 1997), 26--55, Sympos.\ Math., XXXIX, Cambridge Univ.\ Press, Cambridge, 1999.

\bibitem{BB92} M. T. Barlow and R. F. Bass,
\emph{Transition densities for {B}rownian motion of the
{S}ierpinski carpet,}\ Probab.\ Theory Relat.\ Fields, \textbf{91}
(1992), 307--330.


\bibitem{BB3} M. T. Barlow and R. F. Bass,
\emph{Stability of parabolic Harnack inequalities.}\
Trans. Amer. Math. Soc., \textbf{356} (2004), 1501--1533.



\bibitem{BaBaKu1} M.~T.~Barlow, R.~F.~Bass and T.~Kumagai, {\it Stability of parabolic {H}arnack inequalities on
metric measure spaces.} J. Math. Soc. Japan, to appear.

\bibitem{BaKi97} M.~T.~Barlow and J.~Kigami, {\it Localized eigenfunctions of the Laplacian on \pcf\ self-similar
sets,} J.\ London Math.\ Soc., {\bf 56} (1997), no.\ 2, 320--332.

\bibitem{CRS} P. Ciatti, F. Ricci, M. Sundari, {\it Heisenberg-Pauli-Weyl Uncertainty Inequalities
and Polynomial Volume Growth,} to appear in Advances in Mathematics.








\bibitem{FuSh92} M.~Fukushima and T.~Shima, {\it On a spectral analysis for the {S}ierpinski gasket,}
Potential Anal., {\bf 1} (1992), 1--35.



\bibitem{FoSi} G.~B.~Folland and A.~Sitaram, {\it The uncertainty principle: A mathematical survey,} J.\ Fourier
Anal.\ Appl., {\bf 3} (1997), no.\ 3, 207--238.


\bibitem{GriTel01} A.~Grigo{\'r}yan and A.~Telcs,
{\it Sub-{G}aussian estimates of heat kernels on infinite graphs,} Duke math.\ J., {\bf 109} (2001), no.\ 3, 451--510.

\bibitem{GriHuLa03} A.~Grigo{\'r}yan, J.~Hu and K.~Lau, {\it Heat kernels on metric measure spaces and an
application to semilinear elliptic equations,} Trans.\ Amer.\ Math.\ Soc., {\bf 355} (2003), no.\ 5, 2065--2095.

\bibitem{Gri04} A.~Grigo{\'r}yan,  {\it Heat kernel upper bounds on fractal spaces,} preprint.

\bibitem{HaKu1} B.~M.~Hambly and T.~Kumagai, {\it Heat kernel estimates for symmetric random walks on a class of
fractal graphs and stability under rough isometries,} Proc.\
Symp.\ Pure Math., {\bf 72} (2004), 233--259.

\bibitem{Hei01} J.~Heinonen, ``Lectures on analysis on metric spaces,''  Universitext, Springer-Verlag, New York, 2001.

\bibitem
{KL} J. Kigami and M. L. Lapidus,
\emph{Weyl's problem for the spectral distribution of Laplacians on p.c.f.
self-similar fractals.}\
Comm. Math. Phys. {\bf 158} (1993), 93--125.



\bibitem{Ki1} J.  Kigami,
\emph{Harmonic calculus on
\pcf\  self--similar sets.}\  Trans.  Amer.  Math.  Soc.
\textbf{335} (1993),
 721--755.

\bibitem{Ki2} J.  Kigami,
\emph{Effective resistances for harmonic structures on
\pcf\  self--similar sets.}\  Math. Proc. Cambridge Philos. Soc.
\textbf{115} (1994),
291--303.

\bibitem{Ki01} J.~Kigami, ``Analysis on Fractals,'' Cambridge University Press, New York, 2001.

\bibitem
{Ki03} J.~Kigami,
\emph{Harmonic analysis for resistance forms.} J. Functional Analysis {\bf204} (2003), 399--444.


\bibitem{Ki04} J.  Kigami,
\emph{Local Nash inequality and inhomogeneity of heat kernels,}\    Proc. London Math. Soc. (3)  \textbf{89}  (2004),   525--544.

\bibitem{KrTe03} B.~Kr\"{o}n and E.~Teufl, {\it Asymptotics of the transition probabilities of the simple random
walk on self-similar graph,} Trans.\ Amer.\ Math.\ Soc.,
{\bf356} (2003), no.\ 1, 393--414.


\bibitem
{KZ92} S. Kusuoka and X. Y. Zhou,
\emph{Dirichlet forms on fractals: Poincar\'e constant and resistance.}\
Probab. Theory Related Fields \textbf{93} (1992), 169--196.

\bibitem{KZ98} S. Kusuoka and X. Zhou,
\emph{Waves on fractal-like manifolds and effective energy propagation.}\
Probab. Theory Related Fields \textbf{110}
(1998), 473--495.

\bibitem{MT1} L. Malozemov and A. Teplyaev,
\emph{Pure point spectrum of the Laplacians on fractal graphs.}\
J. Funct. Anal.  {\bf129} (1995), 390--405.


\bibitem{MT2} L. Malozemov and A. Teplyaev,
\emph{Self-similarity, operators and dynamics.}\
Math. Phys. Anal. Geom. {\bf6} (2003), 201--218.





\bibitem{OkSt04} K.~A.~Okoudjou and R.~S.~Strichartz, {\it Weak uncertainty principles on fractals,}
J.\ Fourier Anal.\ Appl., {\bf 11} (2005), no.\ 3, 315--331.


\bibitem{PrRa85} J.~F.~Price and P.~C.~Racki, {\it Local uncertainty inequalities for Fourier series,}
Proc.\ Amer.\
Math.\ Soc., {\bf 93} (1985), 245--251.


\bibitem{Sab04a} C.~Sabot, {\it Electrical networks, symplectic reductions, and application to the
renormalization map of self-similar lattices.}
Fractal geometry and applications: a jubilee of Benoit Mandelbrot. Part 1,  155--205,
Proc. Sympos. Pure Math. {\bf 72}, Amer. Math. Soc.,  2004.


\bibitem{Sab04b} C.~Sabot, {\it Laplace operators on fractal lattices with random blow-ups.}
Potential Anal.  {\bf 20} (2004) 177--193.

\bibitem{Sab03} C.~Sabot, {\it Spectral properties of self-similar lattices and iteration of rational maps.}
M\'em. Soc. Math. Fr. (N.S.) No. 92 (2003), vi+104 pp.

\bibitem{Sab00} C.~Sabot, {\it Pure point spectrum for the Laplacian on unbounded nested fractals.}
J. Funct. Anal.  {\bf 173} (2000) 497--524.





\bibitem{S-C02b}  L.~Saloff-Coste,  "Aspects of Sobolev-type inequalities,"  London Mathematical Society
Lecture Note Series, {\bf289}. Cambridge University Press,  2002.

\bibitem{Str89} R.~S.~Strichartz, {\it Uncertainty principles in harmonic analysis,} J.\ Funct.\ Anal., {\bf
84} (1989), 97--114.

\bibitem{Str96} R.~S.~Strichartz, {\it Fractals in large,} Can.\ J.\ Math., {\bf 50} (1996), no.\ 3, 638--657.

\bibitem{Str99} R.~S.~Strichartz, {\it Analysis on fractals,} Notices Amer.\ Math.\ Soc.,
{\bf 46} (1999), 1199--1208.

\bibitem{Str03} R.~S.~Strichartz, {\it Function spaces on fractals,} J.\ Funct.\ Anal., {\bf 198} (2003) 43--83.


\bibitem{Str03f} R.~S.~Strichartz, {\it  Fractafolds based on the Sierpinski gasket and their spectra,}
Trans. Amer. Math. Soc. {\bf355} (2003), 4019--4043.

\bibitem{Str05} R.~S.~Strichartz, {\it Fractal Differential Equations: A Tutorial}, Princeton University Press, 2006.


\bibitem{Tep98} \mbox{A. Teplyaev{,} }\emph{Spectral Analysis on Infinite Sierpi\'nski Gaskets,}\,
J. Funct. Anal., \textbf{159} (1998), 537-567.

\bibitem{vscc}
N.~Th~Varopoulos, L.~Saloff-Coste, T.~Coulhon,
``Analysis and geometry on groups,''
Cambridge Tracts in Mathematics, 100, Cambridge University Press, Cambridge, 1992.
\end{thebibliography}
\end{document}